\documentclass{amsart}
\usepackage {amssymb}
\usepackage {amsmath}
\usepackage {graphicx}

\newtheorem*{prop}{Proposition}

\begin{document}
\title[Algorithm of construction of all knots.]{Algorithm 
of construction of all knots, links with given 
number of crosses on diagram of knot, link, using braids.}
\author[S. S. Serova]{Svetlana S. Serova}
\author[S. A. Serov]{Serge A. Serov}
\address{Sarov, Russia}
\email{saserov@rol.ru}


\begin{abstract}
Algorithm of construction of all knots, 
links with given number of crosses on diagram of knot, link is offered. 
This algorithm is based on simple proposition, 
that there is a representation of knot (link) as closure of braid 
with n threads and length of this braid does not exceed n(4n-5)+2.
\end{abstract}

\maketitle

\section{Introduction}

Classification of knots, here and below we consider knots, links and braids 
in $R^{3}$, is one of the basic problems of the theory of knots, starting 
with works of Tait \cite{Tait898}. To classify knots means somehow or other 
to order a 
set of all knots, for example, to divide a set of all knots on equivalence 
classes to any reference indication, property of knots and to build all 
knots (or to specify really working algorithm of build-up of all knots) for 
the given equivalence class. It seems natural as such reference indication 
to take the least number of crosses of a knot (at an identification of a 
knot with its diagram -- a projection to some plane). By virtue of 
Alexander's known theorem (see, for example, \cite{Prasolov97}) any knot can 
be presented 
as closure of some braid. Braids are organized much easier knots. Therefore, 
it is natural to use braids for classification of knots. Before describing 
of proposed algorithm of construction of knots (links) with the given number 
of crosses on the diagram of a knot (link), we shall make some estimates.

\section{Preliminary estimates}

In computer program, implementing proposed algorithm, it is supposed to use 
partially ``br9z.p'' program in the language Pascal, written by Short and 
Morton in 1985. ``Br9z.p'' program is taken from the site of the Liverpool 
mathematical group \underline{www.liv.ac.uk}. 
In this program, authors used connection of group of knots $B_{n} $ 
with permutation group $S_{n} $. 
The permutation group of a set with $n$ elements consists of $n!$ 
elements. 
Therefore, in ``br9z.p'' program for braids coding the 
two-dimensional array $m$ by a size $\left( {n!} \right) \times n$ of 
numbers from 0 up to $n$ was used. As already $9! = 362880$ is very major 
number, Short and Morton have been forced to be restricted to a maximum 
number of threads of braids $n = 9$ (in ``jones12.p'' program, written in 
1994, Morton increased number of threads up to 12 by usage of a special 
representation of the braid group $B_{n} $ on $2^{n}$ dimensional space). 
Using results of Vogel \cite{Vogel90}, we shall prove, 
that in the problem of 
classification of knots (links) it is possible to substitute 
$4n^{2} - 5n + 2$ for $n!$. For comparison, we may mark, that 
$4 \cdot 9^{2} - 5 \cdot 9 + 2 = 281$ is much less $9! = 362880$.

\begin{prop} There is a representation of a knot (link) as 
closure of braid with $n$ threads and length of this braid (i.e. an 
amount of generators in braid representation as product of standard 
generators $b_{1}, \dots, b_{n - 1} $ of the group 
$B_{n} $) does not exceed $4n^{2} - 5n + 2$, where 
$n$ -- number of crosses for the given diagram of a knot (link).
\end{prop}
\begin{proof}
The diagram of (oriented) knot, link as a result of the 
 {\it operation of erasure of crosses} disintegrates on (oriented) 
closed curves, named  {\it Seifert circles,} - see fig.~\ref{operation} 
(fig.~\ref{operation} is taken from \cite{Prasolov97}). 
Each cross can break up no more than to 2 Seifert circles. 
Therefore, after an erasure of all diagram crosses we shall have 
$s \leqslant 2n$ Seifert circles, where $n$ number of crosses.

\begin{figure}[htbp]
  \begin{center}
    \includegraphics{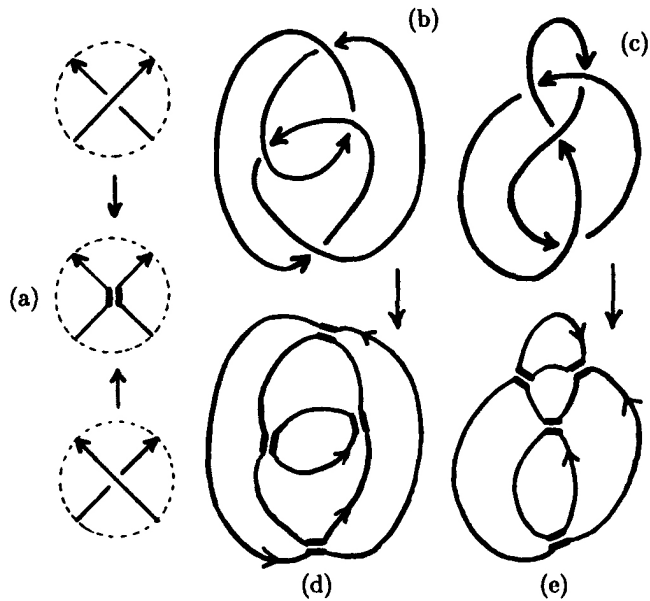}
    \caption{The operation of erasure of crosses.}
    \label{operation}
  \end{center}
\end{figure}

This estimate is strict in the sense, that there are links, for which 
$s = 2n$: for example, $n$-component link, consisting of $n$ eight curves, 
represented on the top or on the bottom of fig.~\ref{link}, after an erasure 
of all crosses, obviously, disintegrates on $2n$ Seifert circles.

\begin{figure}[htbp]
  \begin{center}
    \includegraphics{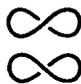}
    \caption{$2$-component link.}
    \label{link}
  \end{center}
\end{figure}

In \cite{Vogel90}, Vogel's article is cited from the book \cite{Prasolov97}, 
Vogel proposed the algorithm of representation of link as closure of braid, 
and proved the estimate: 
obtained braid has $n$ threads and its length (i.e. an amount of 
generators in braid representation as product of standard generators 
$b_{1}, \dots, b_{n - 1} $ of the group $B_{n} $) does not exceed 
$n + \left( {s - 1} \right)\left( {s - 2} \right)$, 
where $n$ -- number of crosses on the diagram of knot (link). 
Assuming in Vogel's estimate $s \leqslant 2n$, we 
can assert, that there is a representation of a knot (link) as closure of 
braid with $n$ threads and braid length does not exceed 
$n + \left( {s - 1} \right)\left( {s - 2} \right) \leqslant 4n^{2} - 5n + 2$, 
where $n$ -- number of crosses of the given diagram of knot (link). 
\end{proof}

Incidentally, it is possible to mark, that proposed by Vogel algorithm of 
representation of a link as closure of braid, from our point of view, is 
difficult for realization in practice for two reasons: 

1) a starting point of algorithm is coding of knot (link) -- very intricate 
problem, actually substituting classification of knots;

2) the {\it operation of replacement of infinity} (see, for example, 
\cite{Prasolov97}), used in the algorithm, is difficult for programming.

If the problem is restricted to classification only knots, then coding of 
representing braids becomes greatly simpler, as it is possible to be 
restricted to braids, which correspond to cyclical permutation of the order 
$n$ -- see above.

\section{Description of algorithm}

The algorithm of construction of all knots (links) with the given number $n$ 
of crosses on the diagram of a knot (link) can be formulated as follows.

1. In the braids group $B_{n} $, $n \geqslant 2$ with standard generators 
$b_{1}, \dots, b_{n - 1} $ we select words with $4n^{2} - 5n + 2$ characters. 
If the problem is restricted to classification only knots, we select words 
with $4n^{2} - 5n + 2$ characters, which corresponds cyclical permutation of 
the order $n$.

2. Using Duhornoy algorithm of reduction \cite{Duhornoy95}, 
we reject superfluous braids -- we select only one representative 
from each equivalence class of braids 
(it is necessary to pay attention, that the algorithm of reduction can 
increase a word length).

3. We use Markov relations (see, for example, \cite{Prasolov97}) for final 
selection of braids, giving at closure different knots (links).

Item 2, apparently, without major detriment, can be omitted. Leaving in 
items 2-3 braids with minimum word length (minimum amount of crosses), as a 
result of closure of selected braids we shall receive diagrams of knots 
(links) with minimum number of crosses. As a result of work of suggested 
algorithm we shall build all knots and links, which diagrams have 
$m \leqslant n$ crosses. Therefore, if we are interested only in knots (links) 
with the least number of crosses on the diagram $m = n$, we should use 
algorithm twice: first to find a set of knots (links) $K_{ \leqslant n - 1} 
$ with the least number of crosses on the diagram $m \leqslant n - 1$, then 
to find a set of knots (links) $K_{ \leqslant n} $ with the least number of 
crosses on the diagram $m \leqslant n$, then the set of knots (links) 
$K_{n} $ with the least number of crosses on the diagram $m = n$ is equal 
to set difference $K_{n} = K_{ \leqslant n} - K_{ \leqslant n - 1} $.

\section*{Conclusion}

In conclusion, it would be desirable to mark, that at preparation this 
article the report on Third Khariton School Readings \cite{Serova03} has been 
used. 
With necessary for understanding of problems discussed above theoretical 
questions on the theory of knots it is possible to familiarize, for example, 
in chapters 1-3 of the book \cite{Prasolov97}.

The program, implementing proposed algorithm, and results of its work is 
supposed to be described in the following article.


\begin{thebibliography}{1}

\bibitem{Tait898} {\bf P. G. Tait},
{\it On Knots I, II, III}, 
Scientific Papers, vol. 1, 
Cambridge University Press, London (1898) 273--347.

\bibitem{Prasolov97} {\bf V. V. Prasolov and A. B. Sosinsky},
{\it Knots, links, braids and three-dimensional manifolds}, 
(in Russian), MCNMO, Moscow (1997).

\bibitem{Vogel90} {\bf P. Vogel},
{\it Representation of links by braids: A new algorithm}, 
Comment. Math. Helvetici. {\bf 65} (1990) 104--113.

\bibitem{Duhornoy95} {\bf P. Duhornoy},
{\it From large cardinals to braids via distributive algebra}, 
J. Knot Theory Ram. {\bf 4} (1995) 33--79.

\bibitem{Serova03} {\bf S. S. Serova},
{\it About geometry and physics of knots, link, braids. 
Classification of knots and links}, 
The Report on Third Khariton School Readings, 
(in Russian, not published), Sarov, Russia (2003).

\end{thebibliography}
\end{document}